\documentclass[10pt,conference,letterpaper]{IEEEtran}

\IEEEoverridecommandlockouts

\usepackage{times,amsmath,epsfig,amssymb,graphicx,amsfonts,amsthm}
\usepackage{subfigure}
\usepackage{empheq}
\usepackage{graphicx}
\usepackage{psfrag}
\usepackage{fge}
\usepackage{amsmath}
\usepackage{amsthm}
\usepackage{color}
\usepackage{amsfonts}   
\usepackage{amssymb}    
\usepackage{mathrsfs}
\usepackage{setspace}
\usepackage{dblfloatfix}
\usepackage{tikz}
\usepackage{tkz-tab}
\usepackage{algorithm}
\usepackage{algcompatible}
\usepackage{lipsum}
\usepackage{color}
\usepackage{mathrsfs}
\usepackage{epstopdf}
\usepackage{enumerate}
\usepackage{cite}

\floatname{algorithm}{Algorithm}
\newtheorem{myfact}{Fact}
\newtheorem{mydef}{Definition}
\newtheorem{myassump}{Assumption}
\newtheorem{mytheorem}{Theorem}

\newtheorem{mycorollary}{Corollary}
\newtheorem{myproposition}{Proposition}

\newtheorem{myobj}{Objective}

\newcounter{ale}

\newenvironment{liste}{\begin{itemize}}{\end{itemize}}
\newcommand{\aliste}{\begin{liste} \setcounter{ale}{1}}
\newcommand{\zliste}{\end{liste}}

\title{{\LARGE {\bf Sensor Placement for Optimal Kalman Filtering: \\ Fundamental Limits, Submodularity, and Algorithms}}}
\author{V.~Tzoumas, A.~Jadbabaie, G.~J.~Pappas{$^{\star}$}
\thanks{$^{\star}$All authors are with the Department of Electrical and Systems Engineering, University of Pennsylvania, Philadelphia, PA 19104-6228 USA (email: {\fontsize{8}{8}\selectfont\ttfamily\upshape \{vtzoumas, jadbabai, pappasg\}@seas.upenn.edu}).}
\thanks{This work was supported in part by TerraSwarm, one of six centers of STARnet, a Semiconductor Research Corporation program sponsored by MARCO and DARPA, in part by AFOSR Complex Networks Program and in part by AFOSR MURI CHASE.}
}

\usepackage{tikz}
\begin{document}
\maketitle

%
%
%
%
%


\begin{abstract}
In this paper, we focus on sensor placement in linear dynamic estimation, where the objective is to place a small number of sensors in a system of interdependent states so to design an estimator with a desired estimation performance.
In particular, we consider a linear time-variant system that is corrupted with process and measurement noise, and study how the selection of its sensors affects the estimation error of the corresponding Kalman filter.  Our contributions are threefold:  First, we prove among other design limits that the number of sensors grows linearly with the system's size for fixed minimum mean square estimation error and number of output measurements over an observation interval.  Second, we prove that the $\log\det$ of the error covariance of the Kalman filter with respect to the system's initial condition and process noise is a supermodular and non-increasing set function in the choice of the sensor set.
Third, we provide efficient approximation algorithms that select a small number sensors so to design an optimal Kalman filter with respect to this estimation error ---the worst-case performance guarantees of these algorithms are provided as well. Finally, we illustrate the efficiency of our algorithms using the problem of monitoring of $CO_2$ sequestration sites studied in \cite{weimer2008relaxation}.
\end{abstract}


\begin{tikzpicture}[overlay, remember picture]
\path (current page.north east) ++(-4,-0.2) node[below left] {
	This paper has been accepted for publication in the American Control Conference.
};
\end{tikzpicture}
\begin{tikzpicture}[overlay, remember picture]
\path (current page.north east) ++(-1.7,-0.6) node[below left] {
	Please cite the paper as:
V.~Tzoumas, A.~Jadbabaie, G.~J.~Pappas,
	``Sensor Placement for Optimal Kalman Filtering: 
};
\end{tikzpicture}
\begin{tikzpicture}[overlay, remember picture]
\path (current page.north east) ++(-2.7,-1) node[below left] {
	Fundamental Limits,
	Submodularity, and Algorithms'', American Control Conference (ACC), 2016.
};
\end{tikzpicture}


\section{Introduction}\label{sec:Intro}

In this paper, we aim to monitor dynamic, interdependent phenomena, that is, phenomena with temporal and spatial correlations ---with the term ``spatial'' we refer to any kind of interdependencies between the phenomena.
For example, the temperature at any point of an indoor environment depends across time ---temporal correlation--- on the temperatures of the adjacent points ---spatial correlation \cite{Madsen199567}.  Therefore, these correlations allow to monitor such phenomena using a reduced number of sensors; this is an important observation when operational constraints, such as limited bandwidth and communication power, necessitate the design of estimators using a small number of sensors \cite{rowaihy2007survey, hero2011sensor}.  Hence, in this paper we consider to place a few sensors so to monitor this kind of phenomena.  To this end, we also account for unknown interdependencies, disturbances and inputs in their dynamics \cite{bertsekas2005dynamic}, and so we consider the presence of process noise, i.e., noise that affects directly these dynamics.  In addition, we account for noisy sensor measurements \cite{kailath2000linear}, and so we consider the presence of measurement noise.

Specifically, we consider phenomena modelled as a linear time-variant system that is corrupted with process and measurement noise, and study how the selection of its sensors affect the minimum mean square error of the corresponding Kalman filter \cite{kailath2000linear}.
To this end, we consider that each of the sensors measures a single state of the system.  Thereby, this study is an important distinction in the sensor placement literature in linear systems \cite{clark2012, lin2013design, matni2014regularization,olshevsky2014minimal, 2014arXiv1401.4209P, bullo2014, summers2014submodularity, 2014arXiv1409.3289T, sergio2015minimal, DBLP:journals/corr/TzoumasRPJ15, barabasi2015spectrum, zhao2015gramian}, since the Kalman filter is the optimal linear estimator ---in the minimum mean square sense--- given a sensor set \cite{kalman1960new}.  In particular, this is the first time that the minimum mean square error of the Kalman filter is studied in this literature, where the objective is to place a small number of sensors in a system so to design an estimator with a desired performance ---this objective is combinatorial and, as a result, also important for large-scale systems, since their size necessitates the design of estimators using a small number of sensors \cite{joshi2009sensor}.  
Our contributions are threefold:

First, we identify fundamental limits in the design of the Kalman filter with respect to its sensors.  In particular, given a fixed number of output measurements over an observation interval, we prove among others that the number of sensors grows linearly with the system's size for fixed minimum mean square error ---this is a design limit, especially for complex systems where the system's size is large.  Moreover, given a fixed number of sensors, we prove that the number of output measurements increases only logarithmically with the system's size for fixed estimation error ---the number of output measurements is proportional to the length of the observation interval, therefore the same result holds for the length of the observation interval as well.  Notwithstanding, it also decreases logarithmically with the number of sensors for fixed system's size and estimation error.  Overall, these results quantify the trade-off between the number of sensors and that of output measurements so to achieve a \mbox{specified value for the estimation~error.}

These results are the first to characterize the effect of the sensor set on the minimum mean square error of the Kalman filter.  In particular, in \cite{bullo2014,olshevsky2015eig}, the authors quantify only the trade-off between the total energy of the consecutive output measurements and the number of its selected sensors; that is, they study no estimator with respect to its sensors.  Similarly, in \cite{barabasi2015spectrum}, the authors consider only the maximum-likelihood estimator for the system's initial condition and only for a special class of stable linear time-invariant systems.  Moreover, they consider systems that are corrupted merely with measurement noise, which is white and Gaussian.  Finally, they also assume an infinite observation interval, that is, infinite number of consecutive output measurements.  On the other hand, we assume a finite observation interval and study the Kalman estimator both for the system's initial condition and for the system's state at the time of the last output measurement.  In addition, we consider general linear time-variant systems that are corrupted with both process and measurement noise, of any distribution (with zero mean and finite variance).  Overall, our results are the first to characterize the effect of the sensor set on the minimum mean square error of the Kalman filter, that is, the optimal linear estimator.

Second, we identify properties for the $\log\det$ of the error covariance of the Kalman filter with respect to the system's initial condition and process noise as a sensor set function ---the design of an optimal Kalman filter with respect to  the system's initial condition and process noise implies the design of an optimal Kalman filter with respect to the system's state.  Specifically, we prove that it is a supermodular and non-increasing set function in the choice of the sensor set.  

In contrast, in \cite{krause2008near}, the authors study sensor placement for monitoring static phenomena with only spatial correlations.  To this end, they prove that the mutual information between the chosen and non-chosen locations is submodular.  On the other hand, we consider dynamic phenomena with both spatial and temporal correlations.  Although the temporal correlation can be translated to a spatial one (cf.~Section \ref{sec:Prelim}), this is a fundamental distinction: the temporal correlation is known, since it is governed by the dynamics of the phenomena.  In particular, this extra knowledge allows for the characterization of a richer class of estimation performance metrics, beyond that of mutual information.  Specifically, in this paper we prove that the $\log\det$ of the error covariance of the Kalman filter ---that is, the optimal linear estimator--- with respect to the phenomena's \mbox{initial condition and process noise is supermodular.}

In the sensor scheduling literature \cite{gupta2006stochastic}, only the $\log\det$ of the error covariance of the Kalman filter with respect to the system's state has been studied: \cite{shamaiah2010greedy} and \cite{jawaid2015submodularity}, the authors prove that it is a supermodular set function for special cases of systems and for zero process noise.  Instead, we consider non-zero process noise, and for any system we prove that the $\log\det$ of the error covariance of the Kalman filter with respect to the initial condition and process noise is supermodular.  In addition, the sensor scheduling literature assumes that a sensor set is already placed, and its objective is at each time step to select a possibly different subset of the system's outputs so to optimize an estimation metric \cite{hero2011sensor}. 
However, in the context of sensor placement no sensor set is assumed already placed and, in particular, the objective is to select one that remains fixed over time and achieves a desired estimation performance.  
In this paper we select a fixed sensor set so to optimize the $\log\det$ of the error covariance of the Kalman filter with respect to the system's initial condition and process noise. We prove it to {be supermodular in the choice of the sensor set.}\footnote{In \cite{jawaid2015submodularity}, the authors prove with a counterexample in the context of sensor scheduling that the minimum mean square error of the Kalman filter with respect to the system's state is not in general a supermodular set function.  We can extend this counterexample in the context of minimal sensor placement as well: the minimum mean square error of the Kalman  with respect to the system's state is not in general a supermodular set function with respect to the choice of the sensor set.}
 
Third, we consider two problems of sensor placement for the design of an optimal Kalman estimator ---we refer to the $\log\det$ of the error covariance of the Kalman filter with respect to the system's initial condition and process noise as \emph{$\log\det$ error}.  First, we consider the problem of designing an estimator that guarantees a specified  $\log\det$ error and uses a minimal number of sensors ---we refer to this problem as \ref{pr:min_set}.
Second, we consider the problem of designing an estimator that uses at most $r$ sensors and  minimizes the  $\log\det$ error ---we refer to this problem as \ref{pr:min}.
Naturally, \ref{pr:min_set} and \ref{pr:min} are combinatorial, and in particular, they involve the minimization of a supermodular set function, that is, the $\log\det$ estimation error of the Kalman filter with respect to the system's initial condition and process noise.
Because the minimization of a general supermodular function is NP-hard \cite{Feige:1998:TLN:285055.285059}, we provide efficient approximation algorithms for their general solution, along with their worst-case performance guarantees.
Specifically, we first provide an efficient algorithm for \ref{pr:min_set} that returns a sensor set that satisfies the estimation guarantee of \ref{pr:min_set} and has cardinality up to a multiplicative factor from the minimum cardinality sensor sets that meet the same estimation bound.  Moreover, this multiplicative factor depends only logarithmically on the problem's \ref{pr:min_set} parameters.  And next, we provide an efficient algorithm for \ref{pr:min}  that returns a sensor set of cardinality $l\geq r$ ---$l$ is chosen by the designer--- and achieves a near optimal value for increasing $l$.  Specifically, for $l=r$, it achieves a worst-case approximation factor $1-1/e$.\footnote{Such algorithms, that involve the minimization of supermodular set functions,
are also used in the machine learning \cite{krause2008beyond}, path planning for information acquisition \cite{singh2009efficient, atanasov2014information}, leader selection \cite{clark2012,clark2014_1, clark2014_2}, sensor scheduling \cite{shamaiah2010greedy,jawaid2015submodularity}, actuator placement \cite{olshevsky2014minimal, summers2014submodularity, 2014arXiv1409.3289T, DBLP:journals/corr/TzoumasRPJ15, DBLP:journals/corr/TzoumasJP15, zhao2015gramian} and sensor placement in static phenomena \cite{krause2008near,2015arXiv150600747J} literature.  Their popularity is due to their simple implementation ---they are greedy algorithms--- and provable worst-case approximation factors, that are the best one can achieve in polynomial time for several classes of functions \cite{nemhauser1978best,Feige:1998:TLN:285055.285059}.}

In comparison, the related literature has focused either a) on the optimization of the $\log\det$ of the error covariance of the Kalman filter with respect to the system's state and only for special cases of systems or for zero process noise, or b) on heuristic algorithms that provide no worst-case performance guarantees, or c) on static phenomena.  In particular, in \cite{joshi2009sensor}, the authors minimize the $\log\det$ of the error covariance of the Kalman filter with respect to the state for the case where there is no process noise in the system's dynamics ---to the contrary, in our framework we assume both process and measurement noise.  Moreover, to this end they use convex relaxation techniques that provide no performance guarantees.  Furthermore, in \cite{dhingra2014admm} and \cite{pfister2014}, the authors design an $H_2$-optimal estimation gain with a small number of non-zero columns.  To this end, they also use convex relaxation techniques  that provide no performance guarantees.  In addition, in \cite{2015arXiv150305968B}, the author designs an output matrix with a desired norm so to minimize the minimum mean square error of the corresponding Kalman estimator.  Nevertheless, the resultant output matrix does not guarantee a small number of selected sensors.  Finally, in \cite{das2008sensor}, the authors consider the problem of sensor placement for monitoring static phenomena with spatial correlation.  To this end, they place a small number of sensors so to minimize a worst-case estimation error of an aggregate function, such as the average.  In contrast, we consider dynamic phenomena with both spatial and temporal correlations that, as mentioned in the preceding paragraphs, offer a richer information setting.  In addition, by minimizing the $\log\det$ error of the Kalman filter with respect to the system's initial condition and process noise, we allow for the efficient estimation of any aggregate function.  Overall, with this paper we are the first to optimize the $\log\det$ error of the Kalman filter in the general dynamic case using a small number of sensors and, at the same time, to provide worst-case performance guarantees.

The remainder of this paper is organized as follows.  In Section \ref{sec:Prelim}, we introduce the system, estimation and sensor placement framework, along with our sensor placement problems.  In Section \ref{sec:limitations}, we provide a series of design and performance limits, and characterize the properties of the Kalman estimator with respect to its sensor set; in Section \ref{sec:sub}, we prove that  the $\log\det$ estimation error of the Kalman filter with respect to the system's initial condition and process noise is a supermodular and non-increasing set function in the choice of the sensor set; and in Section \ref{sec:alg}, we provide  approximation algorithms for selecting a few sensors to design an optimal Kalman filter with respect to its $\log\det$ estimation error ---the worst-case performance guarantees of these algorithms are provided as well.   Finally, 
in Section \ref{sec:sim}, we illustrate our analytical findings, and test the efficiency of our algorithms, using simulation results from an integrator chain network and the problem of surface-based monitoring of $CO_2$ sequestration sites studied in \cite{weimer2008relaxation}.
Section \ref{sec:conc} concludes the paper.

\section{Problem Formulation} \label{sec:Prelim}

\paragraph*{Notation}

We denote the set of natural numbers $\{1,2,\ldots\}$ as $\mathbb{N}$, the set of real numbers as  $\mathbb{R}$,  and the set $\{1, 2, \ldots, n\}$ as $[n]$, where $n \in \mathbb{N}$.  Given a set $\mathcal{X}$, $|\mathcal{X}|$ is its cardinality. 
Matrices are represented by capital letters and vectors by lower-case letters.  For a matrix $A$, $A^\top $ is its transpose and $A_{ij}$ its element located at the $i-$th row and $j-$th column. $\|A\|_2\equiv \sqrt{A^\top A}$ is its spectral norm, and $\lambda_{\min}(A)$ and $\lambda_{\max}(A)$ its minimum and maximum eigenvalues, respectively.  Moreover, if ${A}$ is positive semi-definite or positive definite, we write ${A} \succeq {0}$ and ${A}\succ {0}$, respectively.  
Furthermore, $I$ is the identity matrix ---its dimension is inferred from the context; similarly for the zero matrix $0$.    
Finally, for a random variable $x \in \mathbb{R}^n$, $\mathbb{E}(x)$ is its expected value, and $\mathbb{C}(x)\equiv \mathbb{E}\left(\left[x-\mathbb{E}(x)\right]\left[x-\mathbb{E}(x)\right]^\top \right)$ its covariance.
The rest of our notation is introduced when needed.

\subsection{System and Estimation Framework}\label{subsec:system_estim}

For $k \geq k_0$, consider the linear time-variant system
\begin{equation}\label{eq:dynamics}
\begin{split}
x_{k+1} &= A_kx_k + B_ku_k+w_k,\\
y_k &= C_kx_k+D_ku_k+v_k,
\end{split}
\end{equation}
where $x_k \in \mathbb{R}^n$ ($n \in \mathbb{N}$) is the state vector, $y_k \in \mathbb{R}^c$ ($c \in [n]$) the output vector, $u_k \in \mathbb{R}^p$ ($p \in \mathbb{N}$) the input vector, $w_k$ the process noise and $v_k$ the measurement noise.  The initial condition is $x_0$.
Without loss of generality, $k_0$ and $u_k$ are zero.  

\begin{myassump}[For all $k\geq 0$, the initial condition, the process noise and the measurement noise are uncorrelated random variables]\label{ass:initial_cond}
$x_0$ is a random variable with covariance $\mathbb{C}(x_0)\succ 0$.  
Moreover, for all $k\geq 0$, $\mathbb{C}(w_k)\succ 0$ and $\mathbb{C}(v_k)=\sigma^2I$, where $\sigma > 0$.  Finally, for all $k,k' \geq 0$ such that $k\neq k'$, $x_0$, $w_k$ and $v_k$, as well as, $w_k$, $w_{k'}$, $v_k$ and $v_{k'}$, are uncorrelated.\footnote{This assumption is common in the related literature \cite{joshi2009sensor}, and it translates to a worst-case scenario for the problem we consider in this paper.}
\end{myassump}

For $k \geq 0$, consider the vector of measurements $\bar{y}_k$, the vector of process noises $\bar{w}_k$ and the vector of measurement noises $\bar{v}_k$, defined as follows:  
$\bar{y}_k \equiv (y_{0}^\top , y_{1}^\top , \ldots, y_k^\top )^\top$,
$\bar{w}_k \equiv  (w_{0}^\top , w_{1}^\top , \ldots, w_k^\top )^\top$, and
$\bar{v}_k \equiv  (v_{0}^\top , v_{1}^\top , \ldots, v_k^\top )^\top$;
the vector $\bar{y}_k$ is known, while the $\bar{w}_k$ and $\bar{v}_k$ are not.  

\begin{mydef}[Observation interval and its length]\label{def:obs_int}
The interval $[0,k]\equiv \{0, 1, \ldots, k\}$ is called the \emph{observation interval} of \eqref{eq:dynamics}.  Moreover, $k+1$ is its \emph{length}.
\end{mydef}

Evidently, the length of an observation interval $[0,k]$ equals the number of measurements $y_{0}, y_{1}, \ldots, y_k$.

In this paper, given an observation interval $[0,k]$, we consider the minimum mean square linear estimators for $x_{k'}$, for any $k' \in [0,k]$ \cite{kailath2000linear}.  In particular, \eqref{eq:dynamics} implies
\begin{align}\label{eq:initial_to_output}
\bar{y}_k &= \mathcal{O}_k z_{k-1}+\bar{v}_k,
\end{align}
where $\mathcal{O}_k$ is the $c(k+1)\times n(k+1)$ matrix
$
[ L_0^\top C_0^\top, L_1^\top C_1^\top, \ldots,$ $L_k^\top C_k^\top
]^\top,
$
$L_0$ the $n \times n(k+1)$ matrix $[I, 0]$, $L_i$, for $i\geq 1$, the $n \times n(k+1)$ matrix
$
[A_{i-1}\cdots A_0, A_{i-1}\cdots A_1,$ $ \ldots, A_{i-1}, I, 0],
$
and 
$z_{k-1}\equiv (x_0^\top ,\bar{w}_{k-1}^\top )^\top.
$
As a result, the minimum mean square linear estimate of $z_{k-1}$ is the
$
\hat{z}_{k-1}\equiv \mathbb{E}(z_{k-1})+\mathbb{C}(z_{k-1})$ $\mathcal{O}_k^\top \left(\mathcal{O}_k\mathbb{C}(z_{k-1})\mathcal{O}_k^\top +\sigma^2I\right)^{-1}\left(\bar{y}_k-\mathcal{O}_k\mathbb{E}(z_{k-1})-\mathbb{E}(\bar{v}_k)\right)$; its error covariance is
\begin{align}
\Sigma_{z_{k-1}}&\equiv\mathbb{E}\left(( z_{k-1}-\hat{z}_{k-1})(z_{k-1}-\hat{z}_{k-1})^\top \right)\nonumber\\
&=\mathbb{C}(z_{k-1})-\mathbb{C}(z_{k-1})\mathcal{O}_k^\top \nonumber\\
& \left(\mathcal{O}_k\mathbb{C}(z_{k-1})\mathcal{O}_k^\top +\sigma^2I\right)^{-1}\mathcal{O}_k\mathbb{C}(z_{k-1}) \label{eq:Sigma_z}
\end{align}
and its minimum mean square error
\begin{equation}\label{eq:sq_er_0}
\begin{split}
\text{mmse}(z_{k-1})&\equiv\mathbb{E}\left(( z_{k-1}-\hat{z}_{k-1})^\top (z_{k-1}-\hat{z}_{k-1})\right)\\
&=\text{tr}\left(\Sigma_{z_{k-1}}\right).
\end{split}
\end{equation}
As a result, the corresponding minimum mean square linear estimator of $x_{k'}$, for any $k' \in [0,k]$, is
\begin{equation}\label{eq:est_k}
\hat{x}_{k'}=L_{k'}\hat{z}_{k-1},
\end{equation}
(since $x_{k'}=L_{k'}z_{k-1}$),
with minimum mean square error
\begin{equation}\label{eq:sq_er_1}
\text{mmse}(x_{k'})\equiv\text{tr}\left(\ L_{k'} \Sigma_{z_{k-1}}L_{k'}^\top\right).
\end{equation}
In particular, the recursive implementation of \eqref{eq:est_k} results to the Kalman filtering algorithm \cite{bertsekas2005dynamic}.  

In this paper, in addition to the minimum mean square error of $\hat{x}_{k'}$, we also consider per \eqref{eq:est_k} the estimation error metric that is related to the $\eta$-confidence ellipsoid of $z_{k-1}-\hat{z}_{k-1}$ \cite{joshi2009sensor}.  Specifically, this is the minimum volume ellipsoid that contains $z_{k-1}-\hat{z}_{k-1}$ with probability $\eta$, that is, the $\mathcal{E}_\epsilon\equiv \{z : z^\top \Sigma_{z_{k-1}} z \leq \epsilon\}$,
where $\epsilon \equiv F^{-1}_{\chi^2_{n(k+1)}}(\eta)$ and $F_{\chi^2_{n(k+1)}}$ is the cumulative distribution function of a $\chi$-squared random variable with $n(k+1)$ degrees of freedom \cite{venkatesh2012theory}.  Therefore, the volume of $\mathcal{E}_\epsilon$,
\begin{equation}\label{eq:ellipsoid}
\text{vol}(\mathcal{E}_\epsilon)\equiv \frac{(\epsilon\pi)^{n(k+1)/2}}{\Gamma\left(n(k+1)/2+1\right)}\det\left(\Sigma_{z_{k-1}}^{1/2}\right),
\end{equation}
where $\Gamma(\cdot)$ denotes the Gamma function \cite{venkatesh2012theory}, quantifies the estimation's error of $\hat{z}_{k-1}$, and as a result, for any $k' \in [0,k]$, of $\hat{x}_{k'}$ as well, since per \eqref{eq:est_k} the optimal estimator for $z_{k-1}$ defines the optimal estimator for $x_{k'}$.  

Henceforth, we consider the logarithm of \eqref{eq:ellipsoid},
\begin{equation}\label{eq:vol_ellipsoid}
\log\text{vol}(\mathcal{E}_\epsilon)=\beta+1/2\log\det\left(\Sigma_{z_{k-1}}\right);
\end{equation}
$\beta$ is a constant that depends only on $n(k+1)$ and $\epsilon$, in accordance to \eqref{eq:ellipsoid}, and as a result, we refer to the $\log\det\left(\Sigma_{z_{k-1}}\right)$ as the $\log\det$ estimation error of the Kalman filter of \eqref{eq:dynamics}:

\begin{mydef}[$\log\det$ estimation error of the Kalman filter]
Given an observation interval $[0,k]$, the $\log\det\left(\Sigma_{z_{k-1}}\right)$ is cal- led the \emph{$\log\det$ estimation error of the Kalman filter} of \eqref{eq:dynamics}. 
\end{mydef}

In the following paragraphs, we present our sensor placement framework, that leads to our sensor placement problems. 

\subsection{Sensor Placement Framework}

In this paper, we study the effect of the selected sensors in \eqref{eq:dynamics} on $\text{mmse}(x_0)$ and $\text{mmse}(x_k)$.  Therefore, this translates to the following conditions on $C_k$, for all $k\geq 0$, in accordance with the minimal \mbox{sensor placement literature~\cite{olshevsky2014minimal}.}

\begin{myassump}[$C$ is a full row-rank constant zero-one matrix]\label{assump:Diag_C}
For all $k\geq 0$, $C_k=C \in R^{c \times n}$, where $C$ is a zero-one constant matrix.  Specifically, each row of $C$ has one element equal to one, and each column at most one, such that $C$ has rank $c$.
\end{myassump}

In particular, when for some $i$, $C_{ij}$ is one, the  $j$-th state of $x_{k}$ is measured; otherwise, it is not.
Therefore, the number of non-zero elements of $C$ coincides with the number of placed sensors in \eqref{eq:dynamics}. 

\begin{mydef}[Sensor set and sensor placement]
Consider a $C$ per Assumption \ref{assump:Diag_C} and define $\mathcal{S}\equiv\{i: i\in [n] \text{ and } C_{ji}=1, \text{ for some } j \in [r]\}$; $\mathcal{S}$ is called a \emph{sensor set} or a \emph{sensor placement} and each of its elements a \emph{sensor}.
\end{mydef}


\subsection{Sensor Placement Problems}\label{subsec:problems}

We introduce three objectives, that we use to define the sensor placement problems we consider in this paper.

\begin{myobj}[Fundamental limits in optimal sensor placement]\label{obj:limits}
Given an observation interval $[0,k]$, $i \in \{0,k\}$ and a desired $\emph{\text{mmse}}(x_i)$, identify fundamental limits in the design of the sensor set.
\end{myobj}

As an example of a fundamental limit, we prove that the number of sensors grows linearly with the system's size for fixed estimation error $\text{mmse}(x_i)$ ---this is clearly a major limitation, especially when the system's size is large.  This result, as well as, the rest of our contributions with respect to Objective \ref{obj:limits}, is presented in Section \ref{sec:limitations}.

\begin{myobj}[$\log\det$ estimation error as a sensor set function]\label{obj:sub}
Given an observation interval $[0,k]$, identify properties of the $\log\det\left(\Sigma_{z_{k-1}}\right)$ as a sensor set function.
\end{myobj}

We address this objective in Section \ref{sec:sub}, where we prove that $\log\det\left(\Sigma_{z_{k-1}}\right)$ is a supermodular and non-increasing set function with respect to the choice of the sensor set ---the basic definitions of supermodular set functions are presented in that section as well.

\begin{myobj}[Algorithms for optimal sensor placement]\label{obj:alg}
Given an observation interval $[0,k]$, identify a sensor set $\mathcal{S}$ that solves either the \emph{minimal sensor placement problem:}

\begin{equation}\tag{$\mathcal{P}_1$}\label{pr:min_set}
\begin{aligned}
& \underset{\mathcal{S} \subseteq [n]}{\text{minimize}}  \;
 \quad |\mathcal{S}| \\
&\text{subject to} \quad  \log\det\left(\Sigma_{z_{k-1}}\right) \leq R;
\end{aligned}
\end{equation}
or the \emph{cardinality-constrained sensor placement problem for mi- nimum estimation error:}
\begin{equation}\tag{$\mathcal{P}_2$}\label{pr:min}
\begin{aligned}
& \underset{\mathcal{S} \subseteq [n]}{\text{minimize}} 
 \; \quad \log\det\left(\Sigma_{z_{k-1}}\right) \\
&\text{subject to} 
\quad |\mathcal{S}| \leq r.
\end{aligned}
\end{equation}
\end{myobj}

\ref{pr:min_set}  designs an estimator that guarantees a specified error and uses a minimal number of sensors, and \ref{pr:min} an estimator that uses at most $r$ sensors and minimizes $\log\det\left(\Sigma_{z_{k-1}}\right)$.  The corresponding algorithms are provided in Section \ref{sec:alg}.


\section{Fundamental Limits in Optimal Sensor Placement}\label{sec:limitations}

In this section, we present our contributions with respect to Objective \ref{obj:limits}.  In particular, given $i \in \{0,k\}$ and an observation interval $[0,k]$, we prove that the number of sensors grows linearly with the system's size for fixed $\text{mmse}(x_i)$.  Moreover, given a sensor set of fixed cardinality, we prove that the length of the observational interval increases only logarithmically with the system's size for fixed $\text{mmse}(x_i)$.  Notwithstanding, it also decreases logarithmically with the number of sensors for fixed system's size and $\text{mmse}(x_i)$.  Overall, these novel results quantify the trade-off between the number of sensors and that of output measurements \mbox{so to achieve a specified value for $\text{mmse}(x_i)$.}

To this end, given $i \in \{0,k\}$, we first determine a lower and upper bound for $\text{mmse}(x_i)$.\footnote{The extension of Theorem \ref{th:perfomance_lim} to the case $\mu=1$ is straightforward, yet notationally involved; as a result, we omit it.}

\begin{mytheorem} [A lower and upper bound for the estimation error with respect to the number of  sensors and the length of the observation interval]\label{th:perfomance_lim}
Consider a sensor set $\mathcal{S}$, an observation interval $[0,k]$ and a non-zero $\sigma$.  Moreover, let $\mu\equiv\max_{m\in [0,k]}\|A_m\|_2$ and assume $\mu \neq 1$.  Finally, denote the maximum diagonal element of $\mathbb{C}(x_0)$, $\mathbb{C}(x_0)^{-1}$, and $\mathbb{C}(w_{k'})$, among all $k' \in [0,k]$, as $\sigma_0^{2}$, $\sigma_0^{-2}$, and $\sigma_w^{2}$, respectively.  Given $i \in \{0,k\}$,
\begin{align}
\frac{n\sigma^2l_i}{|\mathcal{S}| \left(1-\mu^{2(k+1)}\right)/\left(1-\mu^2\right)+\sigma^{2}\sigma_0^{-2}} \leq \emph{\text{mmse}}(x_i)\leq nu_i,
\label{eq:lim_0}
\end{align} 
where $l_0=1$, $u_0=\sigma_0^2$, $l_k=\lambda_{\min}\left(L_k^\top L_k\right)$ and $u_k=(k+1)\lambda_{\max}\left(L_k^\top L_k\right)\max(\sigma_0^2,\sigma_w^2)$.
\end{mytheorem}
\begin{proof}
We first prove the lower bound in \eqref{eq:lim_0}: observe first that 
$\text{mmse}(x_0) \geq \text{mmse}(x_0)^{w_{\cdot}=0}$,
where $\text{mmse}(x_0)^{w_{\cdot}=0}$ is the minimum mean square error of $x_{0}$ when the process noise $w_k$ in \eqref{eq:dynamics} is zero for all $k\geq 0$.
To express $\text{mmse}(x_0)^{w_{\cdot}=0}$ in a closed form similar to \eqref{eq:lse_k_final}, note that in this case \eqref{eq:initial_to_output} becomes
$\bar{y}_k = \tilde{\mathcal{O}}_k x_{0}+\bar{v}_k$,
where 
$\tilde{\mathcal{O}}_k \equiv \left[C_0^\top , \Phi_{1}^\top C_{1}^\top , \ldots, \Phi_{k}^\top C_{k}^\top \right]^\top $ and $\Phi_{m}\equiv A_{m-1}\cdots A_{0}$, for $m> 0$, and $\Phi_{m}\equiv I$, for $m=0$.  Thereby, from Corollary E.3.5 of \cite{bertsekas2005dynamic}, the minimum mean square linear estimate of $x_{0}$, denoted as
$
\hat{x}_{k_0}^{w_{\cdot}=0},
$
has error covariance
\begin{align}\label{eq:cov_x_0_state_noise_zero}
\Sigma_{k_0}^{w_{\cdot}=0}&\equiv\mathbb{E}\left(( x_{0}-\hat{x}_{k_0}^{w_{\cdot}=0})(x_{0}-\hat{x}_{k_0}^{w_{\cdot}=0})^\top \right)\nonumber\\
&=\mathbb{C}(x_0)-\mathbb{C}(x_0)\tilde{\mathcal{O}}_k^\top \nonumber\\
& \left(\tilde{\mathcal{O}}_k\mathbb{C}(x_0)\tilde{\mathcal{O}}_k^\top +\sigma^2I\right)^{-1}\tilde{\mathcal{O}}_k\mathbb{C}(x_0) ,
\end{align}
and minimum mean square error
\begin{align}
\text{mmse}(x_0)^{w_{\cdot}=0}&\equiv \text{tr}\left(\Sigma_{k_0}^{w_{\cdot}=0}\right)\nonumber\\
&=\sigma^2\text{tr}\left[\left(\tilde{\mathcal{O}}_k^\top\tilde{\mathcal{O}}_k+\sigma^{2}\mathbb{C}(x_0)^{-1}\right)^{-1}\right]\label{eq:lse_aux_22}\\
&\equiv \sigma^2\text{tr}\left[\left(\tilde{O}_k+\sigma^{2}\mathbb{C}(x_0)^{-1}\right)^{-1}\right]\label{eq:lse_aux_33},
\end{align}
where we deduce \eqref{eq:lse_aux_22} from \eqref{eq:cov_x_0_state_noise_zero} using the Woodbury matrix identity (Corollary 2.8.8 of \cite{bernstein2009matrix}), and \eqref{eq:lse_aux_33} from \eqref{eq:lse_aux_22} using the notation $\tilde{O}_k\equiv\tilde{\mathcal{O}}_k^\top \tilde{\mathcal{O}}_k$.  In particular, $\tilde{O}_k$ is
the observability matrix 
$
\tilde{O}_k = \sum_{m=0}^{k-1}\Phi_{m}^\top C_k^\top C_k\Phi_{m}
$
of \eqref{eq:dynamics} \cite{Chen:1998:LST:521603}.

Hence, $\text{mmse}(x_0) \geq \sigma^2\text{tr}\left[\left(\tilde{O}_k+\sigma^{2}\mathbb{C}(x_0)^{-1}\right)^{-1}\right]$, and since the arithmetic mean of a finite set of positive numbers is at least as large as their harmonic mean, using \eqref{eq:lse_aux_33},
\begin{align*}
\text{mmse}(x_0)\geq \frac{n^2\sigma^2}{\text{tr}\left(\tilde{O}_k+\sigma^{2}\mathbb{C}(x_0)^{-1}\right)} 
\geq\frac{n^2\sigma^2}{\text{tr}\left(\tilde{O}_k\right)+n\sigma^{2}\sigma_0^{-2}}.
\end{align*}
Now, for $i \in [n]$, let $I^{(i)}$ be the $n \times n$ matrix where $I_{ii}$ is one, while $I_{jk}$ is zero, for all $(j,k)\neq (i,i)$. 
Then,
$
\text{tr}(\tilde{O}_k)=\text{tr}\left(\sum_{m=0}^{k}\Phi_{m}^\top C^\top C\Phi_{m}\right)
=\sum_{i=1}^{n} s_i \text{tr}\left(\sum_{m=0}^{k}\Phi_{m}^\top I^{(i)}\Phi_{m}\right)
$;
now,
$
\text{tr}\left(\sum_{m=0}^{k}\Phi_{m}^\top I^{(i)}\Phi_{m}\right)\leq n\lambda_{\max}\left( \sum_{m=0}^{k}\Phi_{m}^\top I^{(i)}\Phi_{m}\right)
= n\|\sum_{m=0}^{k}\Phi_{m}^\top I^{(i)}\Phi_{m}\|_2
\leq n\sum_{m=0}^{k}\|\Phi_{m}\|_2^2,
$
because $\|I^{(i)}\|_2=1$, and from the definition of $\Phi_{m}$ and Proposition 9.6.1 of \cite{bernstein2009matrix},
$
\sum_{m=0}^{k}\|\Phi_{m}\|_2^2
\leq 
\frac{1-\mu^{2(k+1)}}{1-\mu^2}.
$
Therefore, 
$
\text{tr}(\tilde{O}_k)\leq \sum_{i=1}^{n} s_i n\frac{1-\mu^{2(k+1)}}{1-\mu^2} = n|\mathcal{S}| \frac{1-\mu^{2(k+1)}}{1-\mu^2},
$
and as a result, the lower bound in \eqref{eq:lim_0} for $\text{mmse}(x_0)$ follows.

Next, we prove the upper bound in \eqref{eq:lim_0}, using \eqref{eq:lse_aux_333}, which is proved in the proof of Proposition \ref{prop:closed_form_only_output}, and \eqref{eq:sq_er_1} for $k'=0$: $O_k+\sigma^{2}\mathbb{C}(z_{k-1})^{-1} \succeq \sigma^{2}\mathbb{C}(z_{k-1})^{-1}$, and as a result, from Proposition 8.5.5 of \cite{bernstein2009matrix}, $\left(O_k+\sigma^{2}\mathbb{C}(z_{k-1})^{-1}\right)^{-1}\preceq \sigma^{-2}\mathbb{C}(z_{k-1})$.  Hence,
$
\text{mmse}(x_0)\leq\text{tr}\left[L_0\mathbb{C}(z_{k-1})L_0^\top\right]\leq n\sigma_0^2.
$

Finally, to derive the lower and upper bounds for $\text{mmse}(x_k)$, observe that 
$\text{mmse}(x_0)\leq \text{mmse}(z_{k-1})$ and $\text{mmse}(z_{k-1})\leq n(k+1)\max(\sigma_0^2,\sigma_w^2)$
---the proof follows using similar steps as above. Then, from Theorem 1 of \cite{fang1994inequalities},
\begin{align*}
\lambda_{\min}\left(L_k^\top L_k\right)\text{mmse}(z_{k-1})\leq &\text{mmse}(x_k) \leq\\
&  \lambda_{\max}\left(L_k^\top L_k\right)\text{mmse}(z_{k-1}).
\end{align*}
The combination of these inequalities completes the proof.
\end{proof}

The upper bound corresponds to the case where no sensors have been placed ($C=0$). 
On the other hand, the lower bound corresponds to the case where $|\mathcal{S}|$ sensors have been placed.

As expected, the lower bound in \eqref{eq:lim_0} decreases as the number of sensors or the length of the observational interval increases; the increase of either should push the estimation error downwards.   Overall, this lower bound quantifies a fundamental performance limit in the design of the Kalman estimator: for example, given a fixed number of finite sensors and output measurements ---recall that per Definition \ref{def:obs_int} the length of the observation interval equals the number of the output measurements--- the $\text{mmse}(x_0)$ remains bounded away from zero.  In addition, this bound decreases only inversely proportional to the number of sensors, and it increases linearly with the system's size.  For fixed and non-zero $\lambda_{\min}\left(L_k^\top L_k\right)$, this scaling extend to the $\text{mmse}(x_k)$ as well.

\begin{mycorollary}[Trade-off among the number of sensors, estimation error and the length of the observation interval]\label{cor:sensors_tradeoff}
Consider an observation interval $[0,k]$, a non-zero $\sigma$, and that the desired value for $\emph{\text{mmse}}(x_0)$ is $\alpha$ ($\alpha > 0$).  Moreover, let $\mu\equiv\max_{m\in [0,k]}\|A_m\|_2$ and assume $\mu \neq 1$.  Finally, denote the maximum diagonal element of  $\mathbb{C}(x_0)^{-1}$ as $\sigma_0^{-2}$.  Given $i \in \{0,k\}$, any sensor set $\mathcal{S}$ that achieves $\emph{\text{mmse}}(x_i)=\alpha$ satisfies:
\begin{align}\label{eq:limitation_sensors}
|\mathcal{S}| \geq \left(n\sigma^2l_i/\alpha-\sigma^{2}\sigma_0^{-2}\right)\frac{1-\mu^2}{1-\mu^{2(k+1)}},
\end{align}
where $l_0=1$ and $l_k=\lambda_{\min}\left(L_k^\top L_k\right)$.
\end{mycorollary}

That is, the number of necessary sensors increases as the minimum mean square error or the number of output measurements decreases.  Specifically, it decreases exponentially with the number of output measurement.  Notwithstanding, it grows linearly with the system's size for fixed $\text{mmse}(x_0)$.  This is a fundamental design limit, especially when the system's size is large.  For fixed and non-zero $\lambda_{\min}\left(L_k^\top L_k\right)$, the comments of this paragraph extend to the $\text{mmse}(x_k)$ as well ---on the other hand, if $\lambda_{\min}\left(L_k^\top L_k\right)$ varies with the system's size, since $\lambda_{\min}\left(L_k^\top L_k\right) \leq 1$, 
the number of sensors can increase sub-linearly with the system's size for fixed $\text{mmse}(x_k)$.

\begin{mycorollary}[Trade-off among the length of the observation interval, estimation error and the number of sensors]\label{cor:interval_tradeoff}
Consider a sensor set $\mathcal{S}$, a non-zero $\sigma$, and that the desired value for $\emph{\text{mmse}}(x_0)$ is $\alpha$ ($\alpha > 0$).  Moreover, let $\mu\equiv\max_{m\in [0,k]}\|A_m\|_2$ and assume $\mu \neq 1$.   Finally, denote the maximum diagonal element of  $\mathbb{C}(x_0)^{-1}$ as $\sigma_0^{-2}$.   Given $i \in \{0,k\}$, any observation interval $[0,k]$ that achieves $\emph{\text{mmse}}(x_i)=\alpha$ satisfies:
\begin{align}\label{eq:limitation_interval}
k\geq  \frac{\log\left(1-\left(n\sigma^2l_i/\alpha-\sigma^{2}\sigma_0^{-2}\right)(1-\mu^2)/|\mathcal{S}|\right)}{2\log\left(\mu\right)}-1,
\end{align}
where $l_0=1$ and $l_k=\lambda_{\min}\left(L_k^\top L_k\right)$.
\end{mycorollary}

That is, the number of necessary output measurements increases as the minimum mean square error or the number of sensors decreases.  Moreover, in contrast to our comments on Corollary \ref{cor:sensors_tradeoff} and the number of necessary sensors, Corollary \ref{cor:interval_tradeoff} indicates that the number of output measurements increases only logarithmically with the system's size for fixed error and number of sensors.  On the other hand, it also decreases logarithmically with the number of sensors, and this is our final fundamental limit result.

In the following paragraphs, we prove that the $\log\det$ error of the Kalman filter is a supermodular and non-increasing set function in the choice of the sensor set.  We use this result to provide efficient algorithms for the solution of \ref{pr:min_set} and \ref{pr:min}.

\section{Submodularity in Optimal Sensor Placement}\label{sec:sub}

In this section, we present our contributions with respect to Objective \ref{obj:sub}.  In particular, we first derive a closed formula for $\log\det\left(\Sigma_{z_{k-1}}\right)$ and then prove that it is a supermodular and non-increasing set function in the choice of the sensor set ---to the best of our knowledge, this is the first result that relates supermodularity with the $\log\det$ error of the Kalman filter within the framework of sensor placement and for the general case.  This implies that the greedy algorithms for the solution of \ref{pr:min_set} and \ref{pr:min} return efficient approximate solutions \cite{Feige:1998:TLN:285055.285059}; \cite{nemhauser1978best}.  In Section \ref{sec:alg}, we use this supermodularity result and known results from the literature on submodular function maximization \cite{Nemhauser:1988:ICO:42805} to provide efficient algorithms for the solution of \ref{pr:min_set} and \ref{pr:min}.

We now give the definition of a supermodular set function, as well as, that of an non-decreasing set function ---we follow \cite{citeulike:416650} for this material.

Denote as $2^{[n]}$ the power set of $[n]$.

\begin{mydef}[Submodularity and supermodularity]\label{def:sub}
A function $h:2^{[n]}\mapsto \mathbb{R}$ is \emph{submodular} if for any sets $\mathcal{S}$ and $\mathcal{S}'$, with $\mathcal{S} \subseteq \mathcal{S}' \subseteq [n]$, and any $a\notin \mathcal{S}'$,
\[
h(\mathcal{S} \cup \{a\})-h(\mathcal{S})\geq h(\mathcal{S}' \cup \{a\})-h(\mathcal{S}').
\]
A function $h:2^{[n]}\mapsto \mathbb{R}$ is \emph{supermodular} if $(-h)$ is submodular.
\end{mydef}

An alternative definition of a submodular function is based on the notion of non-increasing set functions.

\begin{mydef}[Non-increasing and non-decreasing set function]\label{def:non-inc}
A function $h:2^{[n]}\mapsto \mathbb{R}$ is a \emph{non-increasing set function} if for any $\mathcal{S} \subseteq \mathcal{S}' \subseteq [n]$, $h(\mathcal{S})\geq h(\mathcal{S}')$. Moreover, $h$ is a \emph{non-decreasing set function} if $(-h)$ is a non-increasing set function.
\end{mydef}

Therefore, a function $h:2^{[n]}\mapsto \mathbb{R}$ is submodular if, for any $a\in [n]$, the function $h_a:2^{[n]\setminus\{a\}}\mapsto \mathbb{R}$, defined as $h_a(\mathcal{S})\equiv h(\mathcal{S}\cup \{a\})-h(\mathcal{S})$, is a non-increasing set function.  This property is also called the \emph{diminishing returns property}.

The first major result of this section follows, where we let
\[
O_k\equiv \mathcal{O}_k^\top \mathcal{O}_k,
\]
given an observation interval $[0,k]$.

\begin{myproposition}[Closed formula for the $\log\det$ error as a sensor set function] \label{prop:closed_form_only_output}
Given an observation interval $[0,k]$ and non-zero $\sigma$, irrespective of Assumption \ref{assump:Diag_C},
\begin{align}\label{eq:lse_k_final}
&\log\det\left(\Sigma_{z_{k-1}}\right)=\nonumber\\
& \quad2n(k+1)\log\left(\sigma\right)-\log\det\left(O_k+\sigma^{2}\mathbb{C}(z_{k-1})^{-1}\right).
\end{align}
\end{myproposition}
\begin{proof}
From \eqref{eq:Sigma_z},
\begin{align}
\Sigma_{z_{k-1}}&=\mathbb{C}(z_{k-1})-\mathbb{C}(z_{k-1})\mathcal{O}_k^\top \nonumber\\
& \left(\mathcal{O}_k\mathbb{C}(z_{k-1})\mathcal{O}_k^\top +\sigma^2I\right)^{-1}\mathcal{O}_k\mathbb{C}(z_{k-1})\label{eq:lse_aux_111}\\
&=\left(\sigma^{-2}\mathcal{O}_k^\top \mathcal{O}_k+\mathbb{C}(z_{k-1})^{-1}\right)^{-1}\label{eq:lse_aux_222}\\
&=\sigma^2\left(O_k+\sigma^{2}\mathbb{C}(z_{k-1})^{-1}\right)^{-1}\label{eq:lse_aux_333},
\end{align}
where we deduce \eqref{eq:lse_aux_222} from \eqref{eq:lse_aux_111} using the Woodbury matrix identity (Corollary 2.8.8 of \cite{bernstein2009matrix}, and \eqref{eq:lse_aux_333} from \eqref{eq:lse_aux_222} using the fact that $O_k=\mathcal{O}_k^\top \mathcal{O}_k$. 
\end{proof}

Therefore, the $\log\det\left(\Sigma_{z_{k-1}}\right)$ depends on the sensor set through $O_k$.  Now, the main result of this section follows, where we make explicit the dependence of $O_k$ on the sensor set $\mathcal{S}$.

\begin{mytheorem}[The $\log\det$ error is a supermodular and non-increasing set function with respect to the choice of the sensor set]\label{th:sub}
Given an observation interval $[0,k]$, the
\begin{align*}
&\log\det\left(\Sigma_{z_{k-1}},\mathcal{S}\right)=\\
&\qquad  2n(k+1)\log\left(\sigma\right)-\log\det\left(O_{k,\mathcal{S}}+\sigma^{2}\mathbb{C}(z_{k-1})^{-1}\right): \\
&\quad\qquad\qquad\qquad\qquad\qquad\qquad\qquad\qquad\qquad \mathcal{S} \in  2^{[n]} \mapsto \mathbb{R}
\end{align*}
is a supermodular and non-increasing set function with respect to the choice of the sensor set $\mathcal{S}$.
\end{mytheorem}
\begin{proof}
For $i \in [n]$, let $I^{(i)}$ be the $n \times n$ matrix where $I_{ii}$ is one, while $I_{jk}$ is zero, for all $(j,k)\neq (i,i)$.  Also, let $\bar{\mathbb{C}}\equiv \sigma^{2}\mathbb{C}(z_{k-1})^{-1}$.
Now, to prove that the $\log\det\left(\Sigma_{z_{k-1}},\mathcal{S}\right)$ is non-increasing, observe that
\begin{equation}\label{eq:gramianTOdelta}
O_{k,\mathcal{S}} = \sum_{m=1}^n s_m \sum_{j=0}^{k}L_j^\top I^{(m)}L_j=\sum_{m=1}^n s_m O_{k,\{m\}},
\end{equation} 

Then, for any $\mathcal{S}_1\subseteq\mathcal{S}_2\subseteq [n]$, \eqref{eq:gramianTOdelta} and $O_{k,\{1\}}$, $O_{k,\{2\}}$, $\ldots$, $O_{k,\{n\}} \succeq 0$ imply $O_{k,\mathcal{S}_1} \preceq O_{k,\mathcal{S}_2}$, and as a result, $O_{k,\mathcal{S}_1}+\bar{\mathbb{C}} \preceq O_{k,\mathcal{S}_2}+\bar{\mathbb{C}}$.  Therefore, from Theorem 8.4.9 of \cite{bernstein2009matrix},
$
\log\det\left(O_{k,\mathcal{S}_1}+\bar{\mathbb{C}}\right) \preceq
\log\det\left(O_{k,\mathcal{S}_2}+\bar{\mathbb{C}}\right),
$
and as a result, $\log\det\left(\Sigma_{z_{k-1}}\right)$ is non-increasing. 

Next, to prove that $\log\det\left(\Sigma_{z_{k-1}}\right)$ is a supermodular set function, it suffices to prove that $\log\det\left(O_{k,\mathcal{S}}+\bar{\mathbb{C}}\right)$ is a submodular one.  
In particular, recall that a function $h: 2^{[n]}\mapsto \mathbb{R}$ is submodular if and only if, for any $ a \in [n]$, the function $h_a: 2^{[n]\setminus\{a\}} \mapsto \mathbb{R}$, where $h_a(\mathcal{S})\equiv h(\mathcal{S}\cup \{a\})-h(\mathcal{S})$, is a non-increasing set function.  Therefore, to prove that  $h(\mathcal{S})=\log\det(O_{k,\mathcal{S}}+\bar{\mathbb{C}})$ is submodular, we may prove that the $h_a(\mathcal{S})$ is a non-increasing set function.  To this end, we follow the proof of Theorem 6 in \cite{summers2014submodularity}: first, observe that
\begin{align*}
h_a(\mathcal{S})&= \log\det(O_{k,\mathcal{S}\cup \{a\}}+\bar{\mathbb{C}})-\log\det(O_{k,\mathcal{S}}+\bar{\mathbb{C}})\\
&=\log\det(O_{k,\mathcal{S}}+O_{k,\{a\}}+\bar{\mathbb{C}})-\log\det(O_{k,\mathcal{S}}+\bar{\mathbb{C}}).
\end{align*}
Now, for any $\mathcal{S}_1\subseteq\mathcal{S}_2\subseteq [n]$ and $t \in [0,1]$, define $O(t)\equiv \bar{\mathbb{C}}+O_{k,\mathcal{S}_1}+t(O_{k,\mathcal{S}_2}-O_{k,\mathcal{S}_1})$ and 
$
\bar{h}(t)\equiv \log\det\left(O(t)+O_{k,\{a\}}\right)-\log\det\left(O(t)\right);
$
it is $\bar{h}(0)=h_a(\mathcal{S}_1)$ and $\bar{h}(1)=h_a(\mathcal{S}_2)$.  Moreover, since $d\log\det(O(t)))/dt=\text{tr}\left(O(t)^{-1} dO(t)/dt\right)$ (cf. equation (43) in \cite{petersen2008matrix},
\begin{align*}
\frac{d\bar{h}(t)}{dt}= \text{tr}\left[\left(\left(O(t)+O_{k,\{a\}}\right)^{-1}-O(t)^{-1}\right)O_{k,21}\right],
\end{align*}
where $O_{k,21}\equiv O_{k,\mathcal{S}_2}-O_{k,\mathcal{S}_1}$.  From Proposition 8.5.5 of \cite{bernstein2009matrix}, 
$
\left(O(t)+O_{k,\{a\}}\right)^{-1} \preceq O(t)^{-1},
$
because $O(t)\succ 0$ for all $t \in [0,1]$, since $\bar{\mathbb{C}} \succ 0$, $O_{k,\mathcal{S}_1} \succeq 0$, and $O_{k,\mathcal{S}_2}\succeq O_{k,\mathcal{S}_1}$.  Thereby, from Corollary 8.3.6 of \cite{bernstein2009matrix},
\begin{align*}
\left(\left(O(t)+O_{k,\{a\}}\right)^{-1}-O(t)^{-1}\right)O_{k,21} \preceq 0.
\end{align*}
As a result, ${d\bar{h}(t)}/{dt}\leq 0$, and
\begin{align*}
h_a(\mathcal{S}_2)=\bar{h}(1)=\bar{h}(0)+\int_{0}^{1}\frac{d\bar{h}(t)}{dt}dt\leq \bar{h}(0)=h_a(\mathcal{S}_1).
\end{align*}
Therefore, $h_a(\mathcal{S})$ is a non-increasing set function, and the proof is complete.
\end{proof}

The above theorem states that for any finite observation interval, the $\log\det$ error of the Kalman filter is  a supermodular and non-increasing set function with respect to the choice of the sensor set.  Hence, it exhibits the diminishing returns property: its rate of reduction with respect to newly placed sensors decreases as the already placed sensors increase.  

On the one hand, this property implies a limit in design, since adding new sensors after a first few becomes ineffective for the reduction of the estimation error.  On the other hand, it implies  that greedy approach is effective \cite{nemhauser1978best,Feige:1998:TLN:285055.285059}.  
In the next section, we use the results from the literature on submodular function maximization \cite{Nemhauser:1988:ICO:42805} and provide efficient algorithms for the solution of \ref{pr:min_set} and \ref{pr:min}. 

\section{Algorithms for Optimal Sensor Placement}\label{sec:alg}

In this section, we present our contributions with respect to Objective \ref{obj:alg}.  
\ref{pr:min_set} and \ref{pr:min} are combinatorial, and in particular, we proved in Section \ref{sec:sub} that they involve the minimization of a supermodular set function, that is, the $\log\det$ error.
Therefore, because the minimization of a general supermodular function is NP-hard \cite{Feige:1998:TLN:285055.285059}, in the following paragraphs we provide efficient approximation algorithms for the general solution of \ref{pr:min_set} and \ref{pr:min}, along with their worst-case performance guarantees.  

Specifically, we first provide an efficient algorithm for \ref{pr:min_set} that returns a sensor set that satisfies the estimation bound of \ref{pr:min_set} and has cardinality up to a multiplicative factor from the minimum cardinality sensor sets that meet the same estimation bound.  Moreover, this multiplicative factor depends only logarithmically on the problem's \ref{pr:min_set} parameters.  Next, we provide an efficient  algorithm for \ref{pr:min}  that returns a sensor set of cardinality $l\geq r$ ---$l$ is chosen by the designer--- and achieves a near optimal value for increasing $l$.

To this end, we first present a fact from the supermodular functions minimization literature that we use so to construct an approximation algorithm for \ref{pr:min_set} ---we follow \cite{citeulike:416650} for this material. 
In particular, consider the following problem, which is of similar structure to \ref{pr:min_set}, where $h:2^{[n]}\mapsto \mathbb{R}$ is a supermodular and non-increasing set function:

\begin{equation}\tag{$\mathcal{P}$} \label{pr:X}
\begin{aligned}
& \underset{\mathcal{S} \subseteq [n]}{\text{minimize}}  \;  
 \quad |\mathcal{S}|\\
&\text{subject to} \quad h(\mathcal{S}) \leq R.
\end{aligned}
\end{equation}

The following greedy algorithm has been proposed for its approximate solution, for which, the subsequent fact is true.

\begin{algorithm}
\caption{Approximation Algorithm for \ref{pr:X}.}
\begin{algorithmic}
\REQUIRE $h$, $R$.
\ENSURE Approximate solution for \ref{pr:X}.
\STATE $\mathcal{S}\leftarrow\emptyset$
\WHILE {$h(\mathcal{S}) > R $} \STATE{ 	
    $a_i \leftarrow a'\in \arg\max_{a \in [n]\setminus \mathcal{S}}\left(
	h(\mathcal{S})-h(\mathcal{S}\cup \{a\})
	\right)$\\
	\quad \mbox{} $\mathcal{S} \leftarrow \mathcal{S} \cup \{a_i\}$	
	}
\ENDWHILE
\end{algorithmic} \label{alg:X}
\end{algorithm}

\begin{myfact}\label{fact}
{Denote as $\mathcal{S}^\star$ a solution to \ref{pr:X} and as $\mathcal{S}_0, \mathcal{S}_1, \ldots$ the sequence of sets picked by Algorithm \ref{alg:X}.  Moreover, let $l$ be the smallest index such that $h(\mathcal{S}_l) \leq R$.  Then,
\[
\frac{l}{|\mathcal{S}^\star|} \leq 1+\log \frac{h([n])-h(\emptyset)}{h([n])-h(\mathcal{S}_{l-1})}.
\]}
\end{myfact}

For several classes of submodular functions, this is the best approximation factor one can achieve in polynomial time \cite{Feige:1998:TLN:285055.285059}.  Therefore, we use this result to provide the approximation Algorithm \ref{alg:minimal-leaders} for \ref{pr:min_set}, where we make explicit the dependence of $\log\det\left(\Sigma_{z_{k-1}}\right)$ on the selected sensor set $\mathcal{S}$.
Moreover, its performance is quantified with Theorem \ref{th:minimal}.
\begin{algorithm}
\caption{Approximation Algorithm for \ref{pr:min_set}.}\label{alg:minimal-leaders}
\begin{algorithmic}
\STATE  For $h(\mathcal{S})=\log\det\left(\Sigma_{z_{k-1}},\mathcal{S}\right)$, where $\mathcal{S} \subseteq [n]$, Algorithm \ref{alg:minimal-leaders} is the same as Algorithm \ref{alg:X}.
\end{algorithmic} 
\end{algorithm}

\begin{mytheorem}[A Submodular Set Coverage Optimization for \ref{pr:min_set}]\label{th:minimal}
Denote a solution to \ref{pr:min_set} as $\mathcal{S}^\star$ and the selected set by Algorithm~\ref{alg:minimal-leaders} as $\mathcal{S}$.  Moreover, denote the maximum diagonal element of $\mathbb{C}(x_0)$ and $\mathbb{C}(w_{k'})$, among all $k' \in [0,k]$, as $\sigma_0^{2}$ and $\sigma_w^{2}$, respectively. Then,
\begin{align}
&\log\det\left(\Sigma_{z_{k-1}},\mathcal{S}\right) \leq R,\label{explain:th:minima3}\\
&\frac{|\mathcal{S}|}{|\mathcal{S}^\star|}\leq 1+\log \frac{\log\det\left(\Sigma_{z_{k-1}},\emptyset\right)-\log\det\left(\Sigma_{z_{k-1}},[n]\right)}{R-\log\det\left(\Sigma_{z_{k-1}},[n]\right)}\nonumber\\ &\hspace{7mm}\equiv F_i, \label{explain:th:minimal1}
\end{align}
where $\log\det\left(\Sigma_{z_{k-1}},\emptyset\right)\leq n(k+1)\log\max(\sigma_0^2, \sigma_w^2)$. 
Finally, the computational complexity of Algorithm \ref{alg:minimal-leaders} is $O(n^2(nk)^3)$.
\end{mytheorem}
\begin{proof}
Let $\mathcal{S}_0, \mathcal{S}_1, \ldots$ be the sequence of sets selected by Algorithm~\ref{alg:minimal-leaders} and $l$ the smallest index such that $\log\det\left(\Sigma_{z_{k-1}},\mathcal{S}_l\right)  \leq R$.  Therefore, $\mathcal{S}_l$ is the set that Algorithm \ref{alg:minimal-leaders} returns, and this proves \eqref{explain:th:minima3}.
Moreover, from Fact~\ref{fact},
\begin{align*}
\frac{l}{|\mathcal{S}^\star|}&\leq 
1+\log \frac{\log\det\left(\Sigma_{z_{k-1}},\emptyset\right)-\log\det\left(\Sigma_{z_{k-1}},[n]\right)}{\log\det\left(\Sigma_{z_{k-1}},\mathcal{S}_{l-1}\right)-\log\det\left(\Sigma_{z_{k-1}},[n]\right)}.
\end{align*}
Now, $l$ is the first time that $\log\det\left(\Sigma_{z_{k-1}},\mathcal{S}_l\right) \leq R$, and a result $\log\det\left(\Sigma_{z_{k-1}},\mathcal{S}_{l-1}\right) > R$.  This implies \eqref{explain:th:minimal1}.  

Furthermore,  
$
\log\det\left(\Sigma_{z_{k-1}},\emptyset\right)=\log\det\left(\mathbb{C}_{z_{k-1}}\right),
$
and so from the geometric-arithmetic mean inequality,
\begin{align*}
\log\det\left(\mathbb{C}_{z_{k-1}}\right) &\leq n(k+1)\log \frac{\text{tr}(\mathbb{C}_{z_{k-1}})}{n(k+1)}\\
&\leq n(k+1)\log \frac{n(k+1)\max(\sigma_0^2, \sigma_w^2)}{n(k+1)}\\
&= n(k+1)\log\max(\sigma_0^2, \sigma_w^2).
\end{align*}

Finally, with respect to the computational complexity of Algorithm \ref{alg:minimal-leaders}, note that the \texttt{while} loop is repeated for at most $n$  times.  Moreover, the complexity to compute the determinant of an $n(k+1) \times n(k+1)$ matrix, using Gauss-Jordan elimination decomposition, is $O((nk)^3)$ (this is also the complexity to multiply two such matrices).  Additionally, the determinant of at most $n(k+1)$ matrices must be computed~so
\begin{align*}
\arg\max_{a \in [n]\setminus \mathcal{S}}\left(
    	\log\det\left(\Sigma_{z_{k-1}},\mathcal{S}\right)- \log\det\left(\Sigma_{z_{k-1}},\mathcal{S}\cup \{a\}\right)
    	\right)
\end{align*}
can be computed.  Furthermore, $O(n)$ time is required to find a maximum element between $n$ available.  Therefore, the computational complexity of Algorithm \ref{alg:minimal-leaders} is $O(n^2(nk)^3)$.
\end{proof}

Therefore, Algorithm~\ref{alg:minimal-leaders} returns a sensor set that meets the estimation bound of \ref{pr:min_set}.  The cardinality of this set is up to a multiplicative factor of $F_i$ from the minimum cardinality sensor sets that meet the same estimation bound ---that is, $F_i$ is a worst-case approximation guarantee for Algorithm~\ref{alg:minimal-leaders}.  Additionally, $F_i$ depends only logarithmically on the problem's \ref{pr:min_set} parameters.  Finally, the dependence of $F_i$ on $n$, $R$ and $\max(\sigma_0^2, \sigma_w^2)$ is expected from a design perspective: increasing the network size $n$, requesting a better estimation guarantee by decreasing $R$, or incurring a noise of greater variance, should all push the cardinality of the selected sensor set upwards. 

From a computational perspective, the matrix inversion is the only intensive procedure of Algorithm \ref{alg:minimal-leaders} ---and it is necessary for computing the minimum mean square error of $\hat{x}_i$.  In particular, it requires $O((nk)^3)$ time if we use the Gauss-Jordan elimination decomposition, since $O_k$ in \eqref{eq:lse_k_final} is an $n(k+1) \times n(k+1)$ matrix.  
On the other hand, we can speed up this procedure using the Coppersmith-Winograd algorithm \cite{coppersmith1987matrix}, which requires $O(n^{2.376})$ time for $n \times n$ matrices. Alternatively, we can use numerical methods, which efficiently compute an approximate inverse of a matrix even if its size is of several thousands \cite{Reusken:2001:ADL:587707.587812}.  Moreover, we can speed up Algorithm \ref{alg:minimal-leaders} using the method proposed in \cite{minoux1978accelerated}, which avoids the computation of $\log\det(\Sigma_{z_{k-1}},$ $\mathcal{S})- \log\det(\Sigma_{z_{k-1}},\mathcal{S}\cup \{a\})$ for unnecessary choices of $a$ towards the computation of the 
\begin{align*}
\arg\max_{a \in [n]\setminus \mathcal{S}}\left(
    	\log\det\left(\Sigma_{z_{k-1}},\mathcal{S}\right)- \log\det\left(\Sigma_{z_{k-1}},\mathcal{S}\cup \{a\}\right)
    	\right).
\end{align*}

Next, we develop our approximation algorithm for \ref{pr:min}.  To this end, we first present a relevant fact from the supermodular functions minimization literature ---we follow \cite{citeulike:416650} for this material.
In particular, consider the following problem, which is of similar structure to \ref{pr:min}, where $h:2^{[n]}\mapsto \mathbb{R}$ is a supermodular, non-increasing and non-positive set function:

\begin{equation}\tag{$\mathcal{P}'$} \label{pr:X'}
\begin{aligned}
& \underset{\mathcal{S} \subseteq [n]}{\text{minimize}}  \;
\quad h(\mathcal{S})\\
&\text{subject to} \quad  |\mathcal{S}| \leq r.
\end{aligned}
\end{equation}

Algorithm \ref{alg:X'} has been proposed for its approximate solution, where $l \geq r$, for which, the subsequent fact is true.

\begin{algorithm}
\caption{Approximation Algorithm for \ref{pr:X'}.}
\begin{algorithmic}
\REQUIRE $h$, $l$.
\ENSURE Approximate solution for \ref{pr:X'}.
\STATE $\mathcal{S}\leftarrow\emptyset$, $i \leftarrow 0$
\WHILE {$i < l$} \STATE{ 	
 $a_i \leftarrow a'\in \arg\max_{a \in [n]\setminus \mathcal{S}}\left(
	h(\mathcal{S})-h(\mathcal{S}\cup \{a\})
	\right)$\\
	\quad \mbox{} $\mathcal{S} \leftarrow \mathcal{S} \cup \{a_i\}$, $i \leftarrow i+1$
	}
\ENDWHILE
\end{algorithmic} \label{alg:X'}
\end{algorithm}

\begin{myfact}\label{fact'}
Denote as $\mathcal{S}^\star$ a solution to \ref{pr:X'} and as $\mathcal{S}_0, \mathcal{S}_1, \ldots, \mathcal{S}_l$ the sequence of sets picked by Algorithm \ref{alg:X'}.  Then, for all $l \geq r$,
\[
h(\mathcal{S}_l)\leq \left(1-e^{-l/r}\right)h(\mathcal{S}^\star).
\]
In particular, for $l=r$, $h(\mathcal{S}_l)\leq \left(1-1/e\right)h(\mathcal{S}^\star)$.
\end{myfact}

Thus, Algorithm \ref{alg:X'} constructs a solution of cardinality $l$ instead of $r$ and achieves an approximation factor $1-e^{-l/r}$ instead of $1-1/e$.  For example, for $l=r$, $1-1/e\cong .63$, while for $l=5r$, $1-e^{-l/r}\cong .99$; that is, for $l=5r$, Algorithm \ref{alg:X'} returns an approximate solution that although violates the cardinality constraint $r$, it achieves a value for $h$ that is near to the optimal one.

Moreover, for several classes of submodular functions, this is the best approximation factor one can achieve in polynomial time \cite{nemhauser1978best}.  Therefore, we use this result to provide the approximation Algorithm \ref{alg:min} for \ref{pr:min}, where we make explicit the dependence of $\log\det\left(\Sigma_{z_{k-1}}\right)$ on the selected sensor set $\mathcal{S}$.  Theorem \ref{th:min} quantifies its performance.

\begin{algorithm}[H]
\caption{Approximation Algorithm for \ref{pr:min}.}\label{alg:min}
\begin{algorithmic}
\STATE For $h(\mathcal{S})=\log\det\left(\Sigma_{z_{k-1}},\mathcal{S}\right)$, where $\mathcal{S} \subseteq [n]$, Algorithm \ref{alg:min} is the same as Algorithm \ref{alg:X'}.
\end{algorithmic} 
\end{algorithm}

\begin{mytheorem}[A Submodular Set Coverage Optimization for \ref{pr:min}]\label{th:min}
Denote a solution to \ref{pr:min} as $\mathcal{S}^\star$ and the sequence of sets picked by Algorithm \ref{alg:min} as $\mathcal{S}_0, \mathcal{S}_1, \ldots, \mathcal{S}_l$.  Moreover, denote the maximum diagonal element of $\mathbb{C}(x_0)$ and $\mathbb{C}(w_{k'})$, among all $k' \in [0,k]$, as $\sigma_0^{2}$ and $\sigma_w^{2}$, respectively.  Then, for all $l \geq r$,
\begin{align}\label{ineq:th:min}
\begin{split}
\log\det\left(\Sigma_{z_{k-1}},\mathcal{S}_l\right)\leq(1-&e^{-l/r})\log\det\left(\Sigma_{z_{k-1}},\mathcal{S}^\star\right)+\\
& e^{-l/r}\log\det\left(\Sigma_{z_{k-1}},\emptyset\right),
\end{split} 
\end{align}
where $\log\det\left(\Sigma_{z_{k-1}},\emptyset\right)\leq n(k+1)\log\max(\sigma_0^2, \sigma_w^2)$. 
Finally, the computational complexity of Algorithm \ref{alg:minimal-leaders} is $O(n^2(nk)^3)$.
\end{mytheorem}
\begin{proof}
$\log\det\left(\Sigma_{z_{k-1}},\mathcal{S}\right)-\log\det\left(\Sigma_{z_{k-1}},\emptyset\right)$ is a supermodular, non-increasing and non-positive set function.  Thus, from Fact \ref{fact'} we derive \eqref{ineq:th:min}.

That $\log\det\left(\Sigma_{z_{k-1}},\emptyset\right)\leq n(k+1)\log \max(\sigma_0^2, \sigma_w^2)$, as well as, the computational complexity of Algorithm \ref{alg:min}, follows as in the proof of Theorem \ref{th:minimal}. 
\end{proof}

Algorithm \ref{alg:min} returns a sensor set of cardinality $l\geq r$ and achieves a near optimal value for increasing $l$.  Moreover, the dependence of its approximation level on $n$, $r$ and $\max(\sigma_0^2, \sigma_w^2)$ is expected from a design perspective: increasing the network size $n$, requesting a smaller sensor set by decreasing $r$, or incurring a noise of greater variance should all push the quality of the approximation level downwards. 

Finally, from a computational perspective, our comments on Algorithm \ref{alg:minimal-leaders} apply here as well.

\section{Examples and Discussion}\label{sec:sim}


\subsection{Integrator Chain Network}\label{subsec:integratorChain}

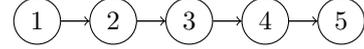
\begin{figure}[t]
\centering
\begin{tikzpicture}
\tikzstyle{every node}=[draw,shape=circle];
\node (v0) at (0:0) {$1$};
\node (v1) at ( 0:1) {$2$};
\node (v2) at ( 0:2) {$3$};
\node (v3) at (0:3) {$4$};
\node (v4) at (0:4) {$5$};

\foreach \from/\to in {v0/v1, v1/v2, v2/v3, v3/v4}
\draw [->] (\from) -- (\to);
\draw
(v0) -- (v1)
(v1) -- (v2)
(v2) -- (v3)
(v3) -- (v4);

\end{tikzpicture}
\caption{A $5$-node integrator chain.}
\label{fig:chain}
\end{figure}

We first illustrate the mechanics and efficiency of Algorithms \ref{alg:minimal-leaders} and \ref{alg:min} using the integrator chain in Fig.~\ref{fig:chain}, where  for all $i \in [5]$, $A_{ii}$ is minus one, $A_{i+1,i}$ one, and the rest of the elements of $A$ are zero.

We run Algorithm \ref{alg:minimal-leaders} for $k \leftarrow 5$, $\mathbb{C}(x_0)$, $\mathbb{C}(w_{k'})$, $\mathbb{C}(v_{k'})\leftarrow I$, for all $k' \in [0,k]$, and $R \leftarrow \log\det(\Sigma_{z_4},\{2,4\})$.  The algorithm returned the sensor set $\{3,5\}$.  Specifically, $\{3,5\}$ is the best such set, as it follows by comparing the values of $\log\det(\Sigma_{z_4},\cdot)$ for all sets of cardinality at most two using MAT- LAB\textsuperscript{\textregistered{}}: for all $i\in [5]$,
$\log\det(\Sigma_{z_4},\{i\})> R$, while for all $i,j \in [5]$ such that $\{i,j\} \neq \{3,5\}$, $\log\det(\Sigma_{z_4},\{i,j\})>\log\det(\Sigma_{z_4},\{3,5\})$.  Therefore, any singleton set does not satisfy the bound $R$,  while $\{3,5\}$ not only satisfies it but also achieves the smallest $\log\det$ error among all other sets of cardinality two; hence, $\{3,5\}$ is the optimal minimal sensor set to achieve the error bound $R$.  Similarly, Algorithm \ref{alg:minimal-leaders} returned the optimal minimal sensor set for every other value of $R$ in the feasible region of \ref{pr:min_set}, $[\log\det(\Sigma_{z_4},[5]),$ $\log\det(\Sigma_{z_4},\emptyset)]$.

We also run Algorithm \ref{alg:min} for $k \leftarrow 5$, $\mathbb{C}(x_0)$, $\mathbb{C}(w_{k'})$, $\mathbb{C}(v_{k'})$ $\leftarrow I$, for all $k' \in [0,k]$, and $r$ being equal to $1, 2,\ldots, 5$, respectively; for all values of $r$ the chosen set coincided with the one of the same size that also minimizes $\log\det(\Sigma_{z_4},\cdot)$; that is, we again observe optimal performance from our algorithms.

Finally, by increasing $r$ from $0$ to $5$, the corresponding minimum value of $\log\det(\Sigma_{z_4},\cdot)$ decreases only linearly from $0$ to approximately $-31$; this is in agreement with the quantitative result of Theorem \ref{th:perfomance_lim}, since for any $M \in \mathbb{R}^{m \times m}$, $\log\det(M) \leq \text{tr}(M)-m$ (Lemma 6, Appendix D in \cite{weimer2010large}), and as a result, for $M$ equal to  $\Sigma_{z_{k-1}}$,
\begin{equation}\label{eq:logdet_trace}
\log\det(\Sigma_{z_{k-1}})\leq \text{mmse}(z_{k-1})-n(k+1),
\end{equation}
while,  for any $k\geq 0$, $ \text{mmse}(x_0)\leq \text{mmse}(z_{k-1})$.

\subsection{$CO_2$ Sequestration Sites}\label{subsec:C02}

We now illustrate the efficiency of Algorithms \ref{alg:minimal-leaders} and \ref{alg:min} using the problem of surface-based monitoring of $CO_2$ sequestration sites \cite{weimer2008relaxation}; these sites are used so to reduce the emissions of $CO_2$ from the power generation plans that burn fossil fuels. In particular, the problem of monitoring $CO_2$ sequestration sites becomes important due to potential leaks.  In addition, since these sites cover large areas, power and communication constrains necessitate the minimal sensor placement for monitoring these leaks, as we undertake in this section.

Specifically, following \cite{weimer2008relaxation}, we consider a) that the sequestration sites form an $9 \times 9$ grid (81 possible sensor locations), and b) the onset of constant unknown leaks. 
Then, for all $k\geq 0$, the $CO_2$ concentration between the sequestration sites is described with the linear time-variant system
\begin{align}\label{eq:C02_sys}
\begin{split}
x_{k+1}&=A_k x_k,\\
y_k&=\left[C, 0\right]x_k+v_k,
\end{split}
\end{align}
where $x_k\equiv (d_k^\top, l_k^\top)^\top$, $d_k \in \mathbb{R}^{81}$ is the vector of the $CO_2$ concentrations in the grid, and $l_k \in \mathbb{R}^{81}$ is the vector of the corresponding leaks rates; $A_k$ describes the dependence of the $CO_2$ concentrations among the adjacent sites in the grid, and with relation to the occurring leaks ---it is constructed as in the Appendix C of \cite{weimer2010large}; and finally, $C \in \mathbb{R}^{81 \times 81}$ is as in Assumption \ref{assump:Diag_C}.  
Hence, next we run Algorithms \ref{alg:minimal-leaders} and \ref{alg:min} so to efficiently estimate the initial condition $x_0$ of \eqref{eq:C02_sys} ---and as a result, detect the constant leaks among the sequestration sites.

\begin{figure}
\centering
\includegraphics[width=\linewidth]{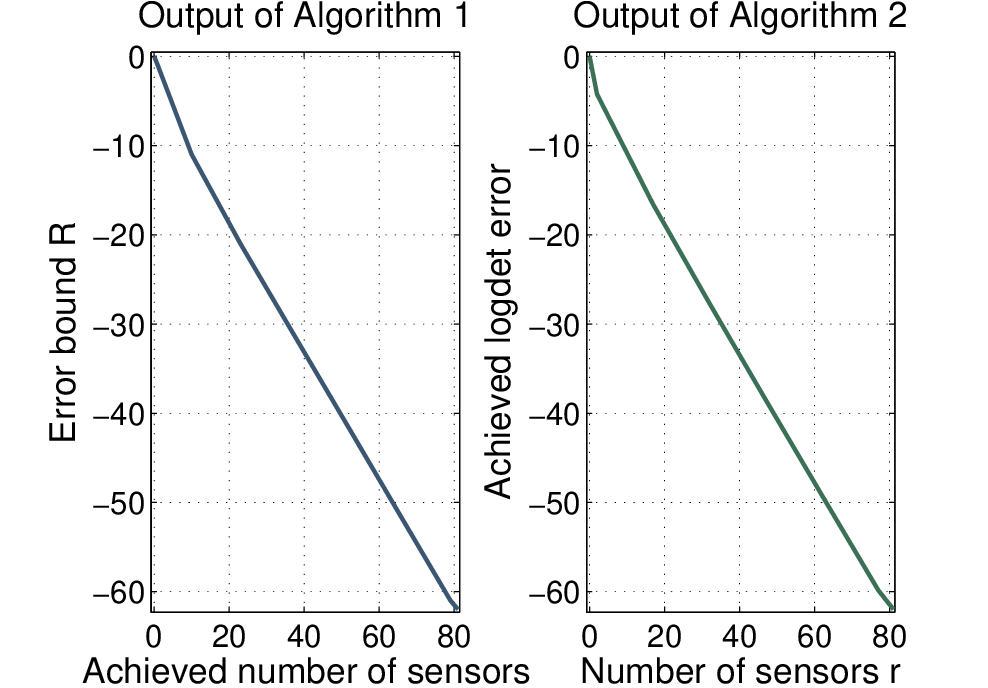}
\caption{Outputs of Algorithms \ref{alg:minimal-leaders} and \ref{alg:min} for the $CO_2$ sequestration sites application of Section \ref{subsec:C02}: In the left plot, the output of Algorithm \ref{alg:minimal-leaders} is depicted for $k \leftarrow 100$, $\mathbb{C}(\bar{x}_0)$, $\mathbb{C}(v_{k'})\leftarrow I$, for all $k' \in [0,k]$, and $R$ that ranged from $\log\det(\Sigma_{z_{99}},[81])\cong -62$ to $\log\det(\Sigma_{z_{99}},\emptyset)=0$ with step size ten.  In the right plot, the output of Algorithm \ref{alg:min} is depicted for $k \leftarrow 100$, $\mathbb{C}(\bar{x}_0)$, $\mathbb{C}(v_{k'})\leftarrow I$, for all $k' \in [0,k]$, as well, and $r$ that ranged from $0$ to $81$ with step size ten.}
\label{fig:Alg1and2together}
\end{figure}

In particular, we run Algorithm \ref{alg:minimal-leaders} for $k \leftarrow 100$, $\mathbb{C}({x}_0)$, $\mathbb{C}(v_{k'})\leftarrow I$, for all $k' \in [0,k]$, and $R$ that ranged from $\log$ $\det (\Sigma_{z_{99}},[81])\cong -62$ to $\log\det(\Sigma_{z_{99}},\emptyset)=0$ with step size ten ---since in \eqref{eq:C02_sys} the process noise is zero, $z_{99}=x_0$.  The corresponding number of sensors that Algorithm \ref{alg:minimal-leaders} achieved with respect to $R$ is shown in the left plot of Fig.~\ref{fig:Alg1and2together}: as $R$ increases the number of sensors decreases, as one would expect when the estimation error bound of \ref{pr:min_set} is relaxed.  

For the same values for $k$, $\mathbb{C}({x}_0)$ and $\mathbb{C}(v_{k'})$, for all $k' \in [0,k]$, we also run Algorithm \ref{alg:min} for $r$ that ranged from $0$ to $81$ with step size ten.  The corresponding achieved values for Problem \ref{pr:min} with respect to $r$ are found in the right plot of Fig.~\ref{fig:Alg1and2together}: as the number of available sensors $r$ increases the minimum achieved value also decreases, as expected by the monotonicity and supermodularity of $\log\det(\Sigma_{z_{99}},\cdot)$ ---since in \eqref{eq:C02_sys} the process noise is zero, $z_{99}=x_0$.  At the same plot, we compare these values with the minimums achieved over a random sample of $80,000$ sensor sets for the various $r$ ---$10,000$ distinct sets for each $r$, where $r$ ranged from $0$ to $81$ with step size ten. Evidently, the achieved minimum random values coincide with those of Algorithm \ref{alg:min}.

In addition, to compare the outputs of Algorithms \ref{alg:minimal-leaders} and \ref{alg:min}, for error bounds $R$ larger than minus twenty, Algorithm \ref{alg:minimal-leaders} returns a sensor set that solves \ref{pr:min_set}, yet does not minimize the $\log\det(\Sigma_{z_{99}},\cdot)$; notwithstanding, the difference in the achieved value with respect to the corresponding output of Algorithm \ref{alg:min} is small.  Moreover, for error bounds  $R$ less than minus twenty, Algorithm \ref{alg:minimal-leaders} returns a sensor set that not only satisfies the bound $R$; it also minimizes $\log\det(\Sigma_{z_{99}},\cdot)$.  Overall, with respect also to our comments in the preceding paragraph, both Algorithms \ref{alg:minimal-leaders} and \ref{alg:min} outperform the theoretical guarantees of Theorems \ref{th:minimal} and \ref{th:min}, respectively.

Finally, as in the example of Section \ref{subsec:integratorChain},  the minimum value of $\log\det(\Sigma_{z_{99}},\cdot)$ ---right plot of Fig.~\ref{fig:Alg1and2together}--- decreases only linearly from $0$ to approximately $-62$ as the number of sensors $r$ increases from $0$ to $81$; this is in agreement with the quantitative result of Theorem \ref{th:perfomance_lim}, in relation to \eqref{eq:logdet_trace} and the fact that $z_{99}=x_0$ ---since in \eqref{eq:C02_sys} the process noise is zero.  Both the integrator chain example of Section \ref{subsec:integratorChain} and the application example of this section exemplify this fundamental design limit presented in Theorem \ref{th:perfomance_lim}: the estimation error of the optimal linear estimator ---the Kalman filter--- decreases only linearly as the number of sensors~increases.

\section{Concluding Remarks}\label{sec:conc}

We considered a linear time-variant system and studied the properties of its Kalman estimator given an observation interval $[0,k]$ and a sensor set $\mathcal{S}$.
Our contributions were threefold.
First, in Section \ref{sec:limitations} we presented several design and performance limits.  For example, we proved that the cardinality of the selected sensors $\mathcal{S}$ grows linearly with the system's size for fixed minimum mean square estimation error and $k$.
Second,  in Section \ref{sec:sub} we proved that the $\log\det$ estimation error of the system's initial condition and process noise is a supermodular and non-increasing set function with respect to the choice of the sensor set for any $k$.  Then, third, in Section \ref{sec:alg}, we used this result to provide efficient approximation algorithms for the solution of \ref{pr:min_set} and \ref{pr:min}, along with their worst-case performance guarantees.  For example, for \ref{pr:min_set}, we provided an efficient algorithm that returns a sensor set that has cardinality up to a multiplicative factor from that of the corresponding optimal solutions; moreover, this factor depends only logarithmically on the problem's parameters.  And for \ref{pr:min}, we provided an efficient algorithm that returns a sensor set of cardinality $l\geq r$ and achieves a near optimal value for increasing $l$.
Finally,  in Section \ref{sec:sim}, we illustrated our analytical findings, and tested the efficiency of our algorithms, using simulation results from an integrator chain network and the problem of surface-based monitoring of $CO_2$ sequestration sites studied in \cite{weimer2008relaxation} ---another application that fits the context of minimal sensor placement for effective monitoring is that of thermal control of indoor environments, such as large offices and buildings \cite{FOldewurtel,sturzenegger}.
Our future work is focused on extending the results of this paper to the problem of minimal sensor placement for sensor scheduling. 

%

\bibliographystyle{IEEEtran}
\bibliography{references}

\end{document}